\documentclass{amsart}
\usepackage{graphicx}
\usepackage{algorithm}
\usepackage{algpseudocode}
\usepackage{booktabs}
\usepackage{multirow}
\usepackage{arydshln}
\usepackage[table,xcdraw]{xcolor}
\usepackage{tabu}
\usepackage{amssymb}
\usepackage{mathtools}
\usepackage{mathrsfs}
\usepackage{amsthm}
\numberwithin{equation}{section}
\numberwithin{figure}{section}
\numberwithin{table}{section}
\usepackage{bbm}
\usepackage{subfig}
\usepackage{enumerate}
\usepackage[section]{placeins}

\usepackage{ifpdf}
\ifpdf
\DeclareGraphicsExtensions{.pdf,.eps,.jpg,.png}	
\usepackage[suffix=]{epstopdf}
\fi
\usepackage{xcolor}
\usepackage[utf8]{inputenc}
\usepackage{hyperref}
\hypersetup{hidelinks}

\long\def\MSC#1\EndMSC{\def\arg{#1}\ifx\arg\empty\relax\else
	{\narrower\noindent%
		{2020 Mathematics Subject Classification}: #1\\} \fi}

\long\def\KEY#1\EndKEY{\def\arg{#1}\ifx\arg\empty\relax\else
	{\narrower\noindent%
		Keywords: #1\\}\fi}

\theoremstyle{plain}

\theoremstyle{definition}

\theoremstyle{remark}

\begin{document}
	\title[tiDAs]{tiDAS: a time invariant approximation of the Delay and Sum algorithm for biomedical ultrasound PSF reconstructions}
	
	
	
	
	\author[C.~Razzetta]{Chiara Razzetta}
	\address[C.~Razzetta]{IRCCS Ospedale Policlinico San Martino, Largo Rosanna Benzi 10, IT-16132 Genova, Italy}
	\email{chiara.razzetta@hsanmartino.it (Corresponding Author)}
	
	\author[S.~Garbarino]{Sara Garbarino}
	\address[S.~Garbarino]{Dipartimento di Matematica, University of Genova, via Dodecaneso 35, IT-16146 Genova, Italy.}
	\email{sara.garbarino@unige.it}
	
	\author[M.~Piana]{Michele Piana}
	\address[M.~Piana]{Dipartimento di Matematica, University of Genova, via Dodecaneso 35, IT-16146 Genova, Italy. IRCCS Ospedale Policlinico San Martino, Largo Rosanna Benzi 10, IT-16132 Genova, Italy}
	\email{michele.piana@unige.it}
	
	\author[M.~Crocco]{Marco Crocco}
	\address[M.~Crocco]{Esaote S.p.A., Via E. Melen 77, IT-16152 Genova, Italy}
	\email{marco.crocco@esaote.com}
	
	\author[F.~Benvenuto]{Federico Benvenuto}
	\address[F.~Benvenuto]{Dipartimento di Matematica, University of Genova, via Dodecaneso 35, IT-16146 Genova, Italy.}
	\email{federico.benvenuto@unige.it}

\maketitle

\begin{abstract}
Ultrasound imaging is a real-time diagnostic modality that reconstructs acoustic signals into visual representations of internal body structures. 
A key component in this process is beamforming, with the Delay and Sum (DAS) algorithm being a standard due to its balance between simplicity and effectiveness. 
However, the computational cost of DAS can be a limiting factor, especially in real-time scenarios where fast frame reconstruction is essential. 
In this work, we introduce a time-invariant approximation of the DAS algorithm (tiDAS), designed to accelerate the reconstruction process without compromising image quality.
By adopting a one-dimensional, row-wise convolutional formulation, tiDAS significantly reduces computational complexity while preserving the core properties of the original model. 
This approach not only enables faster image reconstruction but also provides a structured foundation for the application of deconvolution methods aimed at enhancing resolution. 
Synthetic experiments demonstrate that tiDAS achieves a favorable trade-off between speed and accuracy, making it a promising tool for improving the efficiency of real-time ultrasound imaging.
\end{abstract}


\section{Introduction}
Medical ultrasound imaging is a critical diagnostic tool widely used across various medical disciplines, including cardiology, obstetrics, and radiology \cite{Cobbold2007, Hoskins2010, Szabo2004}. 
This non-invasive technique exploits high-frequency sound waves to produce real-time images of internal anatomical structures, offering valuable insights into both morphology and function. 
Its portability, safety, and cost-effectiveness make it an indispensable modality in clinical practice, particularly for point-of-care and emergency settings.

At the core of ultrasound imaging lies the process of image reconstruction, which converts raw acoustic signals collected by transducer arrays into interpretable visual representations. 
The accuracy and efficiency of this reconstruction process are essential for producing clinically meaningful images, particularly in real-time applications where diagnostic decisions are made on the spot.
Reconstruction involves a sequence of complex signal processing techniques aimed at enhancing image quality, resolution, and diagnostic reliability. 
Notable approaches include deconvolution methods \cite{Basarab2023, Basarab2022, Zhao2016}, regularization-based frameworks \cite{Besson2018, Basarab2016,bertero2006inverse}, plane-wave imaging techniques \cite{Couture2012, Hashemseresht2022}, and strategies employing spatial coherence and super-resolution principles \cite{Christensen2020, Long2022}.
Among these, beamforming, and specifically the Delay and Sum (DAS) algorithm, remains foundational, particularly due to its compatibility with industrial implementations. 

DAS works by introducing controlled delays to signals received by elements of the transducer array and summing them coherently, focusing the acoustic energy in specific directions. 
This process improves spatial resolution and signal-to-noise ratio (SNR), thereby enhancing the visibility of anatomical features. 
Its conceptual simplicity and effectiveness have made DAS the industry standard in many commercial ultrasound systems, but despite its widespread adoption, the DAS algorithm presents a trade-off between accuracy and computational efficiency. 
Indeed, the ability to steer and focus beams dynamically comes at the cost of high computational load, which can become a bottleneck in applications requiring high frame rates or low-latency feedback. 
Furthermore, the choice of transmission and reception parameters significantly influences the overall image quality \cite{Perrot2021}, adding to the complexity of achieving optimal performance in real-time settings.
Importantly, the extent and sophistication of post-processing that can be applied to ultrasound images is inherently constrained by the speed of the initial reconstruction. 
A faster DAS implementation would enable quicker visual feedback, which is crucial in time-sensitive clinical environments such as emergency diagnostics or intraoperative imaging. 

Beyond reducing latency, accelerating the reconstruction stage also would free up computational resources and time that could be redirected toward advanced post-processing operations. 
These include techniques such as speckle reduction, contrast enhancement, tissue characterization, Doppler analysis, elastography, and super-resolution reconstruction, many of which are critical for enabling more detailed and reliable diagnostic insights \cite{Oglat2018,Park2014, Postema2011,  Sigrist2017}. 
Since these methods are often computationally intensive, an efficient reconstruction phase is a prerequisite for applying them effectively without disrupting the real-time constraints of the imaging workflow. 
Conversely, if the reconstruction step itself is computationally heavy, the available time budget for post-processing shrinks, thereby limiting the application of these advanced techniques and reducing the system’s responsiveness.

In light of these challenges, in this study we propose a time-invariant approximation of DAS (tiDAS) designed to reduce computational complexity while maintaining high image fidelity. 
The method is formulated in a one-dimensional framework, which enables a convolutional representation of the reconstruction process along each image row. 
This structure not only accelerates per-frame processing but also paves the way for the integration of deconvolution techniques aimed at further improving resolution and clarity.
By leveraging the principle of time-invariance, the tiDAS approach preserves the essential focusing behavior of standard DAS while simplifying its computational implementation. 
The resulting speed-up in processing is particularly beneficial for real-time imaging applications, where both high frame rates and diagnostic accuracy are critical. 
Additionally, the structured nature of the model facilitates future extensions toward more advanced reconstruction and post-processing pipelines, thereby enhancing the overall utility of ultrasound as a real-time diagnostic tool.

The paper is structured as follows. In Section \ref{sec:math}, we introduce the mathematical background of beamforming and formally derive the tiDAS model. 
Section \ref{sec:method} presents the synthetic simulation setup, along with a quantitative and qualitative evaluation of the proposed method. 
Final remarks and potential future directions are offered in the concluding section.

\section{Local time invariant approximation of DAS}\label{sec:math}
In ultrasound imaging, the goal is to estimate the tissue reflectivity function, which quantifies the body ability to reflect ultrasound signals. The probe records the pressure wave scattered back by the body internal structures.
On the return path, the reflected ultrasound propagates from the field to the probe, striking the piezoelectric elements.
The active elements record a time-dependent signal that has been attenuated by the body. This signal shares the same impulse response function as in the transmission phase and is proportional to the tissue reflectivity function.
In this context, the body is assumed to have a medium density, and the pressure waves are considered to propagate at a medium velocity. As the pressure wave crosses biological tissues, slight variations in density and wave velocity occur. The tissue reflectivity function captures these local variations at each point in the field.

Let $N : = \left[-\overline{N}, \cdots, \overline{N}-1\right]$ denote the set of active elements, symmetric with respect to the center of the probe.
The reflectivity at a given field point $\vec{r}$ is denoted by $g^{\dagger}(\vec{r})$. We denote by $S^F(\vec{r}, \cdot)$ the focused transmitted signal that passes through point $\vec{r}$  over time.
The pressure wave scattered back from $\vec{r}$ reaches all active elements, resulting in a collection of recorded signals:
\begin{equation}
	S_n(\vec{r}, \cdot) = g^{\dagger}(\vec{r}) S^F(\vec{r}, \cdot) \ast H_n(\vec{r}, \cdot) \quad \forall n \in N,
	\label{eqn:SRX}
\end{equation}
where $H_n(\vec{r}, \cdot)$ is the attenuated impulse response for the $n$-th element.
In general, the field does not contain a single scatterer, but a collection of them, each reflecting part of the transmitted signal. Under Born approximation \cite{brignone2007use}, second-order effects due to mutual interactions among scattered waves are neglected \cite{Nigam1959}.
Therefore, B-mode images are typically reconstructed using the Delay and Sum (DAS) algorithm \cite{McKeighen1977}. The core idea is that signals reflected from a single point reach different probe elements at different times. DAS dynamically focuses the image by applying specific delay profiles for each reconstruction point. Time-of-flight differences are used to align and sum signal samples across the probe, yielding an estimate of the tissue reflectivity function at each point, i.e.,
\begin{equation}
	g(\vec{r}) = \sum_{n \in N} S_n \ast \delta_{D_n({\vec{r}})},
	\label{eqn:appF}
\end{equation}
where $D_n({\vec{r}})$  represents the delay applied to the signal received by element $n$.
In the following, we analyze how these delays influence the point spread function (PSF).
\subsection{Delays Properties}

For the sake of simplicity and without loss of generality, we restrict our analysis to the case of a linear probe, where the transmission is not steered. This choice allows us to describe the behavior of the delays in a simplified yet representative setting. Specifically, we assume that the transmission beam is focused at a given depth and centered along the axis of the probe. The receiving elements are identified by their geometric centers $\vec{x}_n = (x_n, 0)$, and the imaging plane is defined as the vertical plane $y = 0$, so that each target point has coordinates $\vec{r} = (x, z)$.
\begin{figure}
	\centering
	\includegraphics[width=\linewidth]{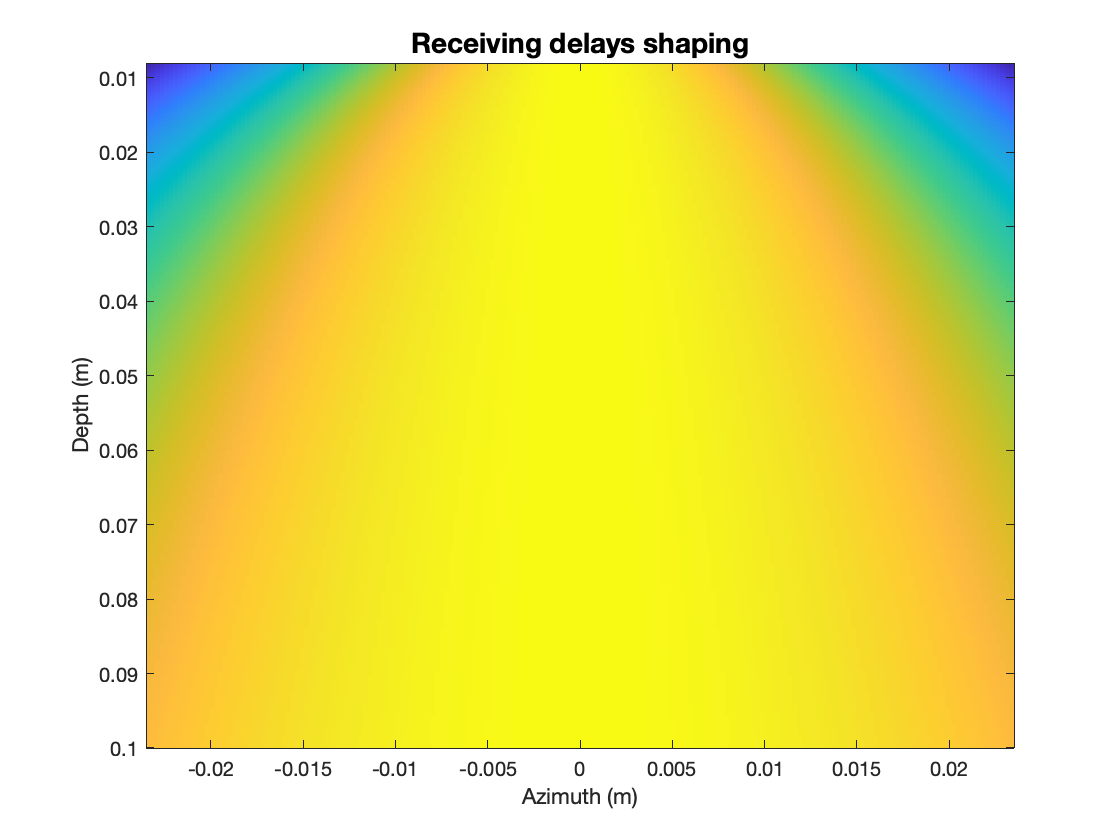}
	\caption{Picturing of all the delays curve obtained for a linear probe of 192 elements making the depth vary between $1 cm $ and $10 cm$. The image points out how the variation is smooth even discarding the fact that usually nearby the probe not all the elements are active.}
	\label{fig:ImageDelays}
\end{figure} 
\begin{figure}
	\centering
	\includegraphics[width=\linewidth]{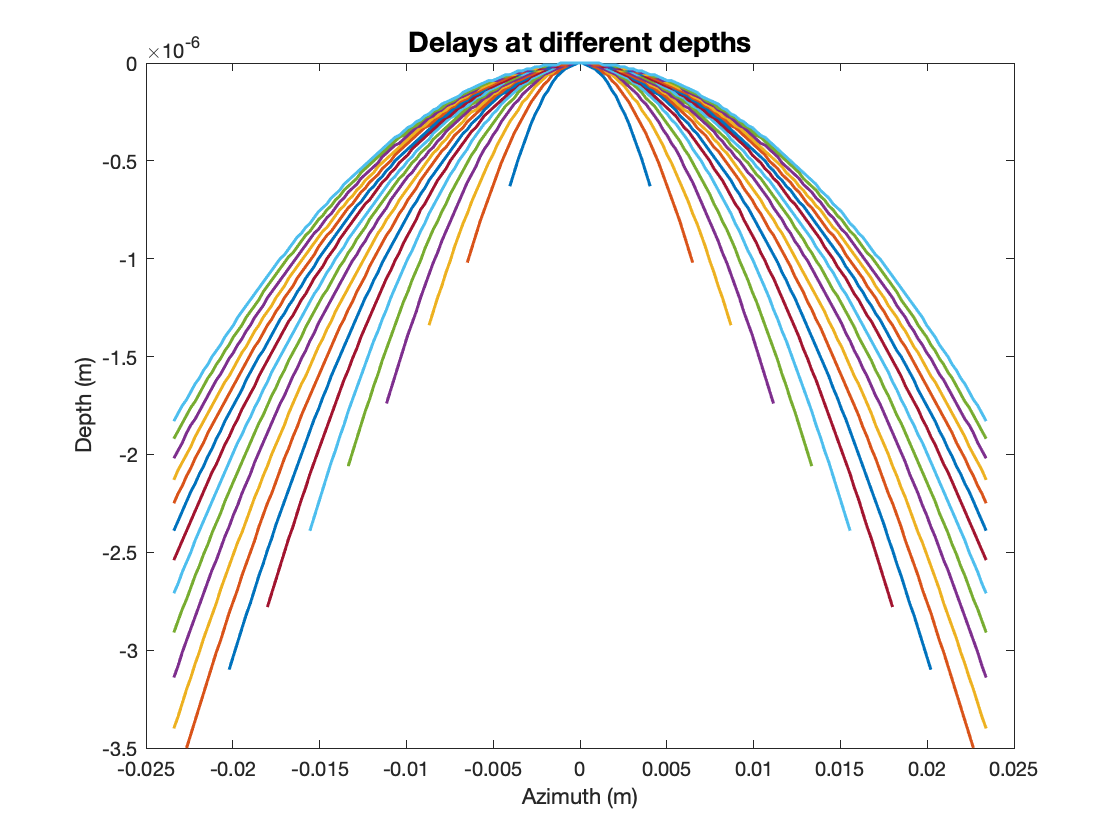}
	\caption{We display different delays curves used to reconstruct points at different depths. For each depth, we take in account also the different aperture estimated by the focal number rule, with focal number set to $1$. The farther the target point is from the probe surface the higher is the flattening of the curve.}
	\label{fig:delAperture}
\end{figure}
In this setting, the receive delay applied to signals from element $n$ when focusing at depth $z$ along the central line $(x=0)$ is given by:
\begin{equation}
    D_n(\vec{r}) \coloneqq \frac{1}{c} \left(z - \sqrt{x_n^2 + z^2} \right),
    \label{eqn:delay}
\end{equation}
where $c$ is the speed of sound in the medium.
This expression reveals that each depth $z$ corresponds to a distinct delay profile across the array, and these delay curves are approximately parabolic in shape. Importantly, these curves vary smoothly with respect to $z$. As shown in Figures~\ref{fig:ImageDelays} and \ref{fig:delAperture}, increasing the imaging depth results in flatter delay curves, meaning the variation across elements becomes less pronounced. Conversely, for shallow depths, fewer elements are active, and the delay curves exhibit more curvature and variation.
To quantify this behavior, we consider the first-order Taylor expansion of the delay function $D_n(\vec{r})$ around a reference depth $\vec{r}_0$:
\begin{equation}
    D_n(\vec{r}) \approx D_n(\vec{r}_0) + \frac{\partial D_n}{\partial z}\Big|_{z_0} (z - z_0),
    \label{eqn:taylor}
\end{equation}
yielding the relative variation:
\begin{equation}
    \frac{D_n(\vec{r}) - D_n(\vec{r}_0)}{D_n(\vec{r}_0)} \approx -\frac{1}{\sqrt{x_n^2 + z_0^2}}(z - z_0).
    \label{eqn:relTaylor}
\end{equation}
This shows that the relative variation of delays decreases:
\begin{enumerate}
    \item As the depth $z_0$ increases, due to increased flattening of the delay curves.
    \item Along the tails of the curve (i.e., for large $|x_n|$), where the denominator grows.
\end{enumerate}
\begin{figure}
	\centering
	\includegraphics[width=\linewidth]{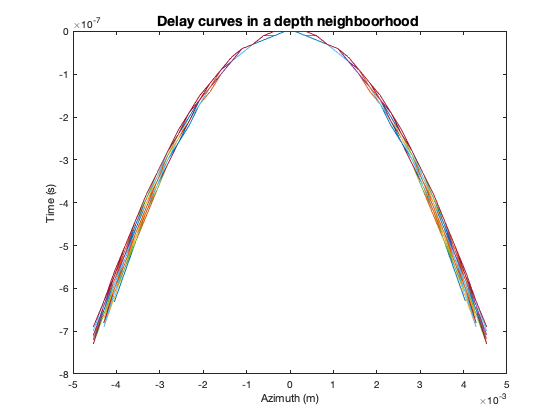}
	\caption{The graph displays many delay curves in the neighborhood of a target point. It is clear nearby the vertex the differences are smaller then along the tails.}
	\label{fig:Neigdelays}
\end{figure}
In Figure~\ref{fig:Neigdelays}, we display a collection of delay curves corresponding to target points within a small neighborhood of a fixed depth. The difference among these curves is minimal near the curve vertex (i.e., for central elements) and more visible along the tails. Notably, these differences become increasingly negligible at greater depths.
These observations support a key assumption in our proposed model: in a small depth neighborhood, the delay curves can be considered approximately invariant. 
This allows us to approximate the DAS algorithm as a time-invariant convolution by fixing the delay profile to that of the transmission focus point and absorbing the residual differences into a multiplicative scaling factor.
The practical implications of this are twofold:
\begin{itemize}
    \item The PSF becomes approximately shift-invariant over small depth intervals, enabling a convolutional model.
    \item The number of required delay profiles to reconstruct an entire image line is drastically reduced, since a single delay profile can be reused within each local neighborhood.
\end{itemize}
In the following subsections, we  formally derive the model underlying the time-invariant approximation of the DAS algorithm.

\subsection{Single scatterer case}
We now focus on the case of a single scatterer to describe the formulation of the PSF. Specifically, we consider a single scatterer located at a point $\vec{r}$, which reflects the incoming ultrasound wave.
The signal received by the probe element $n$ is given by
\begin{equation}
	S_n^{\vec{r}}(t) \coloneqq \left(H^{\vec{r}}_n \ast S^{F, \vec{r}}\right)(t).
	\label{eqn:PSFrx}
\end{equation}
In the DAS aproach, the corresponding reconstructed value of the reflectivity at $\vec{r}$ is:
\begin{equation}
	g^{\vec{r}} = \sum_{n \in N} S_n^{F, \vec{r}} \ast \delta_{D_n^{\vec{r}}}.
	\label{eq:PSF}
\end{equation}
Although equation \eqref{eq:PSF} is written as a convolution, it is important to note that the resulting PSF depends explicitly on time, as the delays $D_n^{\vec{r}}$ are time-dependent.

\subsection{Time invariant DAS (tiDAS)}\label{subsec:tiDAS}
In this section, we introduce an approximation of the standard Delay and Sum (DAS) algorithm, where the receiving delays are assumed to be constant. This leads to a time-invariant version of the DAS, referred to as tiDAS.
To account for the simplification, we incorporate a scaling factor into the PSF, applied in the vicinity of the peak time. This factor acts as a multiplicative weight that depends on the position of the scatterer, thereby adjusting the PSF to better approximate the original time-varying formulation.

To ensure that the weight is only applied over a limited time window, we introduce an indicator function $\chi^{\vec{r}}(t) = \mathbf{1}_{B}(t)$ where the interval $B = [z/c - \epsilon, z/c +\epsilon]$ defines a neighborhood around the expected echo arrival time. Here, $z$ is the axial coordinate of the scatterer, $c$ is the speed of sound in the medium, and $\epsilon$ controls the window width.
The received signal is then modulated by this indicator function, restricting the influence of the scaling to a narrow temporal region of interest, i.e.,
\begin{equation}
	\tilde{g}^{\vec{r}} = \theta_{\vec{r}}\sum_{n \in N} \tilde{S}_n \ast \delta_{D_n},
	\label{eqn:tDAS}
\end{equation}
where $\theta_{\vec{r}} \in \mathbb{R}$ is the multiplicative weight, $\tilde{S}_n = S_n^{F, \vec{r}} \cdot \chi^{\vec{r}}(t)$ is the collection of signals multiplied by the indicator function and $D_n$ is the delay curve corresponding to a point $\vec{s}$ in the domain.
Our goal is to determine the weight $\theta_{\vec{r}}$ and the curve $D_n$ such that
\begin{equation}
	\tilde{g}^{\vec{r}} \approx g^{\vec{r}}. 
\end{equation}
To do so, we compute the norm of the difference between the two reconstructions as
\begin{align}
	\left|\left|\tilde{g}^{\vec{r}} - g^{\vec{r}}\right|\right|_2^2 & = \left|\left|\theta_{\vec{r}}\sum_{n \in N} \tilde{S}_n \ast \delta_{D_n} - \sum_{n \in N} \tilde{S}_n \ast \delta_{D_n^{\vec{r}}}\right|\right|_2^2\nonumber\\
	&=\left|\left|\sum_{n \in N} \tilde{S}_n \ast \left( \theta_{\vec{r}}\delta_{D_n} - \delta_{D_n^{\vec{r}}}\right)\right|\right|_2^2\nonumber.
\end{align}
Applying Plancherel theorem leads to the equality with the same equation in time-frequency domain, i.e.,
\begin{equation}
\begin{split}
    	\left|\left|\tilde{g}^{\vec{r}} - g^{\vec{r}}\right|\right|_2^2 &=\\ &\left|\left|\sum_{n \in N} \hat{\tilde{S}}_n(f) \theta_{\vec{r}} e^{-i 2 \pi f D_n} - \hat{\tilde{S}}_n(f) e^{-i2\pi f D_n^{\vec{r}}}\right|\right|_2^2.
\end{split}
	\label{eqn:diff}
\end{equation}
From now on, we denote by:
$$\begin{cases}
	S_n^f : = \hat{S}_n^{F, \vec{r}}(f)\\
	x_n : = e^{-i 2 \pi f D_n}\\
	y_n : = e^{-i2\pi f D_n^{\vec{r}}}
\end{cases},$$
and:
$$\begin{cases}
    a(f) = \sum_{n \in N} S_n^f x_n\\
    b(f) = \sum_{n \in N} S_n^f y_n
\end{cases}.$$
Thus, it holds:
\begin{align}
    \left|\left|\tilde{g}^{\vec{r}} - g^{\vec{r}}\right|\right|_2^2 &= \left|\left| \theta_{\vec{r}}a(f) - b(f)  \right|\right|_2^2 \nonumber\\&= \int \overline{\left(\theta_{\vec{r}}a(f) - b(f)\right)} \left(\theta_{\vec{r}}a(f) - b(f)\right) \, df \nonumber\\
    &= \theta_{\vec{r}}^2 \left|\left|a(f) \right|\right|_2^2 + \left|\left|b(f) \right|\right|_2^2 +\nonumber\\&- \theta_{\vec{r}} \int \left[\overline{a(f)} b(f) + a(f) \overline{b(f)}\right] \, df
\label{eq:norma}
\end{align}

Deriving equation \ref{eq:norma} by $\theta_{\vec{r}}$ leads to
\begin{align}
    \frac{d}{d \theta_{\vec{r}}}\left|\left|\tilde{g}^{\vec{r}} - g^{\vec{r}}\right|\right|_2^2 &=\nonumber\\  &=2 \theta_{\vec{r}} \left|\left|a(f) \right|\right|_2^2 -  \int \overline{a(f)} b(f) + a(f) \overline{b(f)} \, df ,\nonumber
\end{align}
so that
\begin{align}
    \theta_{\vec{r}} &= \frac{ \int \overline{a(f)} b(f) + a(f) \overline{b(f)} \, df  }{2\left|\left|a(f) \right|\right|_2^2} \nonumber\\
    &= \frac{\int \sum_{n\in N} \left(\overline{S_n^f x_n}\cdot S_n^f y_n + S_n^f x_n \cdot \overline{S_n^f y_n}\right) \, df}{\int \sum_{n\in N} \left(\overline{S_n^f x_n}\cdot S_n^f x_n \right) \, df} \ .
    \label{eqn:theta}
\end{align}

It is worth noticing that we have found a scaling parameter dependent on the received signal, which makes approximation \eqref{eqn:tDAS}, adaptive with respect it.

\subsection{Convolutional model}
Since the approximation depends on the choice of a delay profile and the scale parameter, we can define a convolutional model to reconstruct a whole line of pixels.
Under the tiDAS approximation, the contributions of each transducer element can be aligned using a fixed delay set $D\coloneqq \{ D_n ~:~ n\in N \}$ and scaled by a corresponding set of weights $\Theta$.
This leads to the following convolutional formulation for the reconstruction of a line of pixels:
\begin{equation}
    \mathcal{L} \coloneqq \Theta \sum_{n \in N} \tilde{S}_n \ast \delta_{D_n},
    \label{eqn:line_conv}
\end{equation}
where
\begin{equation}
    \Theta \coloneqq \left\{\theta_{\vec{r}} : \vec{r} \in \mathcal{L} \right\}
\end{equation}
for a fixed delay set $D$.
This formulation becomes feasible due to the definition of the modified reflectivity signal $\tilde{g}^{\vec{r}}$, which removes the explicit dependence on the spatial location $\vec{r}$, thereby allowing the reconstruction to be interpreted as a classical convolution.

By shifting the complexity from time-varying alignment to a convolution with fixed parameters, this approach offers the potential for substantial computational savings over traditional DAS methods, while still preserving the key structural features of the beamformed signal.
The primary challenge, however, lies in the selection of an appropriate delay profile $D_n$, which must effectively account for the geometry of all target points along the reconstruction line. Additionally, the scaling parameters $\Theta$ must be carefully chosen to compensate for spatially varying signal attenuation and focusing effects.

\section{Results}\label{sec:method}
In this section, we present a set of results comparing reconstructions obtained using the standard DAS algorithm and the proposed tiDAS approximation. Our primary objective is to assess the reconstruction accuracy for individual scatterers, with a particular focus on point reflectors.
Next, we analyze the performance of various delay curves and their corresponding scaling parameter sets, evaluating how well they capture the spatial and temporal features of the beamformed signal.
Finally, we demonstrate the effectiveness of the tiDAS model in reconstructing a whole image line, highlighting its potential to reduce computational complexity while maintaining image quality.

All simulations are conducted in Python, using a dedicated simulation framework built on our customized ultrasound imaging software 'parUST' \cite{Razzetta23B, parust}.

\subsection{Simulation setting}
We simulate the behavior of a linear ultrasound probe to compare the performance of the proposed tiDAS approximation with that of the traditional DAS algorithm. The analysis focuses on selected points within the field of view that are geometrically aligned with the central axis of the probe. This choice relies on the assumption of translational invariance along the $x$-axis, which simplifies the computational model without loss of generality for the scenario under study.
The axial direction ($z$-axis) is discretized into $600$ evenly spaced points spanning from $2 \, mm$ to $42\, mm$. This resolution is sufficient to capture fine spatial variations in reconstruction performance.
The number of active elements used during transmission is determined according to Kossoff’s criterion \cite{Kossoff1979}, with a parameter value of $0.6$, which optimizes beam formation by controlling the effective aperture. As for reception, the number of active elements is selected to maintain a fixed focal number of $1$, balancing lateral resolution and focal depth.

To evaluate the imaging response, a single point among the 600 along the z-axis is designated as the reference focal point for the transmission phase, ensuring consistency across all experimental conditions.
In addition to spatial selection, a fixed temporal window is also defined for the analysis.
Specifically, the Full Width at Half Maximum (FWHM) of the PSF at the focal point is used as the time interval, in accordance with the theoretical framework described in Section \ref{sec:math}, where this interval is defined for the correction procedure.
\subsection{FWHM as a natural correction neighborhood}

As introduced in Section~\ref{sec:math}, the tiDAS approximation replaces the dynamically varying receive delays of standard DAS with a fixed delay curve, and introduces a spatially dependent scaling factor applied over a limited time interval. 
To ensure that this correction focuses on the most relevant portion of the signal—namely, the main lobe of the point spread function (PSF)—we propose to define the correction window $B$ using the FWHM of the DAS PSF at the focal point.

The FWHM captures the width of the main lobe where most of the signal energy is concentrated. Using it to define the indicator function $\chi^{\vec{r}}(t)$ ensures that the scaling factor $\theta_{\vec{r}}$ is applied only in a time window that corresponds to the meaningful part of the signal. 
This approach avoids introducing distortion in the side lobes or noise floor and maintains consistency with the underlying physics of beamforming.

To validate this choice, we perform a quantitative evaluation of tiDAS reconstructions for single scatterers across a wide range of depths and acquisition settings. 
Specifically, we focus on two metrics:
\begin{itemize}
    \item The absolute difference in peak amplitude between DAS and tiDAS reconstructions.
    \item The difference in FWHM between the two signals, capturing the similarity in their temporal support.
\end{itemize}

\begin{figure}
	\centering
	\includegraphics[width=\linewidth]{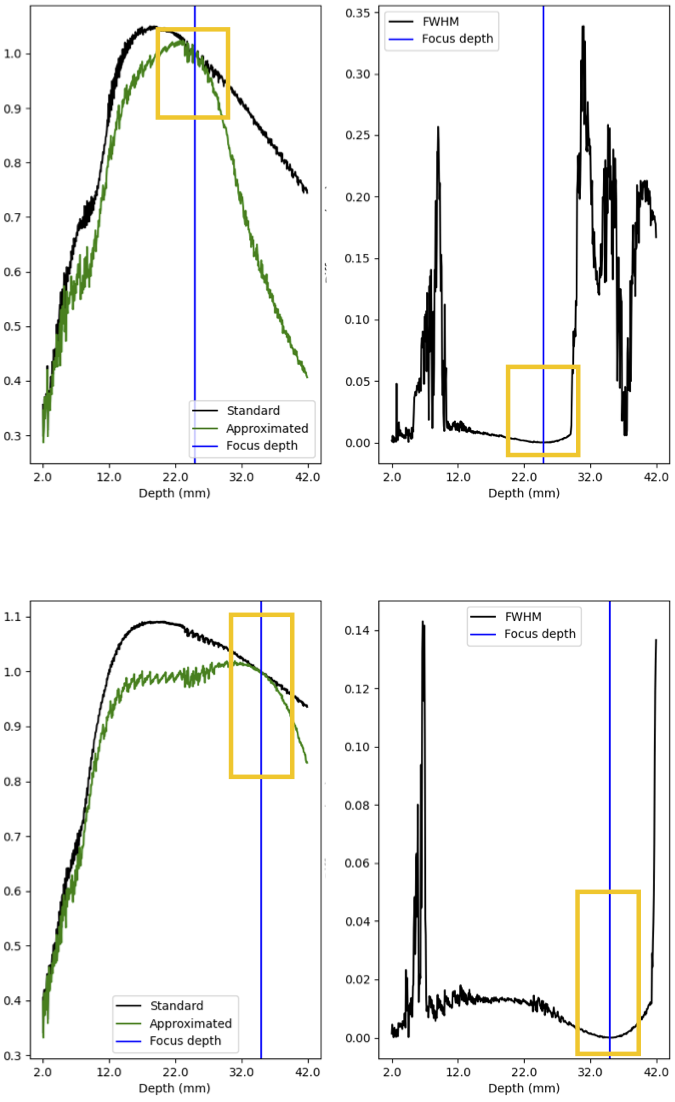}
	\caption{
Panels in the first column show the peak heights of the reconstructions obtained with DAS and tiDAS, normalized with respect to DAS reconstruction in the focal point. 
The peak height of the DAS reconstructions over the depth is represented by the black line, The same curve for the tiDAS algorithm is shown in green. The blue line represents the focal point in transmission.
Panels in the second column show the absolute difference of the FWHM of the DAS and tiDAS reconstructions as a function of depth.
The yellow rectangles highlight the regions where the tiDAS method closely approximates the DAS reconstructions.
The plots in the first and second rows were computed using central frequencies of $7 \, MHz$ and $5 \, MHz$ respectively.}
	\label{fig:k6_fwhm}
\end{figure}

Figure~\ref{fig:k6_fwhm} 
report the ratio between PSF peaks at different depths 

the results for a variety of configurations, including different center frequencies

different active surface of transmission, and focal depths. Each row corresponds to a different focal point and displays both metrics across all 600 scatterer positions. The vertical blue line indicates the focal depth, while the yellow box identifies a neighborhood where the tiDAS approximation remains highly accurate.

\begin{figure}
    \centering
    \includegraphics[width=\linewidth]{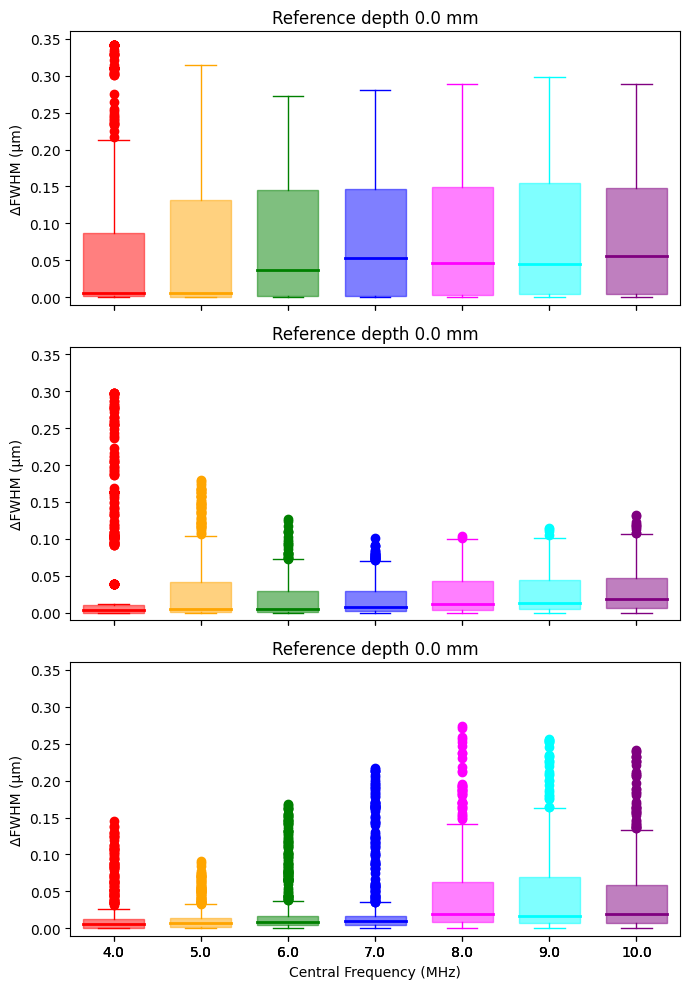}
    \caption{Distribution of FWHM differences between DAS and tiDAS for three focal depths and various center frequencies. Results are aggregated over 600 simulated scatterer positions per configuration.}
    \label{fig:summary_fwhm}
\end{figure}

To consolidate these findings across multiple configurations, we report in Figure~\ref{fig:summary_fwhm} the distribution of FWHM differences between DAS and tiDAS reconstructions, computed over 600 scatterer positions for each central frequency. 
The results are grouped by transmission focus depth—$15$, $25$, and $35\, mm$—and illustrate the effect of both frequency and depth on the accuracy of the approximation.

Each box plot summarizes the distribution of absolute FWHM differences for a given center frequency, with whiskers extending to the 1.5 interquartile range and outliers marked as individual points. 
At $25$ and $35\, mm$, the majority of differences fall below $0.1\, \mu m$, indicating high fidelity in the main lobe shape reconstruction.

These results reinforce the reliability of using the FWHM-defined window as a principled choice for applying the scaling correction $\theta_{\vec{r}}$. 
The window effectively isolates the region of interest around the signal peak, ensuring that the time-invariant approximation captures the essential features of the standard DAS PSF, even under varied imaging conditions.

In the next subsection, we apply this correction scheme across a range of scatterer positions and evaluate the accuracy and stability of the resulting reconstructions.

\subsection{Single points reconstruction}
For each one of the 600 potential scatterer positions within the homogeneous medium, we simulate the presence of a single point scatterer. The corresponding backscattered signals are then reconstructed using both the traditional DAS algorithm and the proposed tiDAS approximation.
This process yields a comprehensive dataset consisting of $600^2$
pairs of reconstructed signals, covering a broad range of spatial configurations. Each signal pair is associated with a corresponding scaling parameter $\theta_{\vec{r}}$,  along with a computed reconstruction error. This dataset provides a robust foundation for quantitatively assessing the performance of tiDAS relative to standard DAS across varying depths and focal positions.
\begin{figure}
    \centering
    \includegraphics[width=\linewidth]{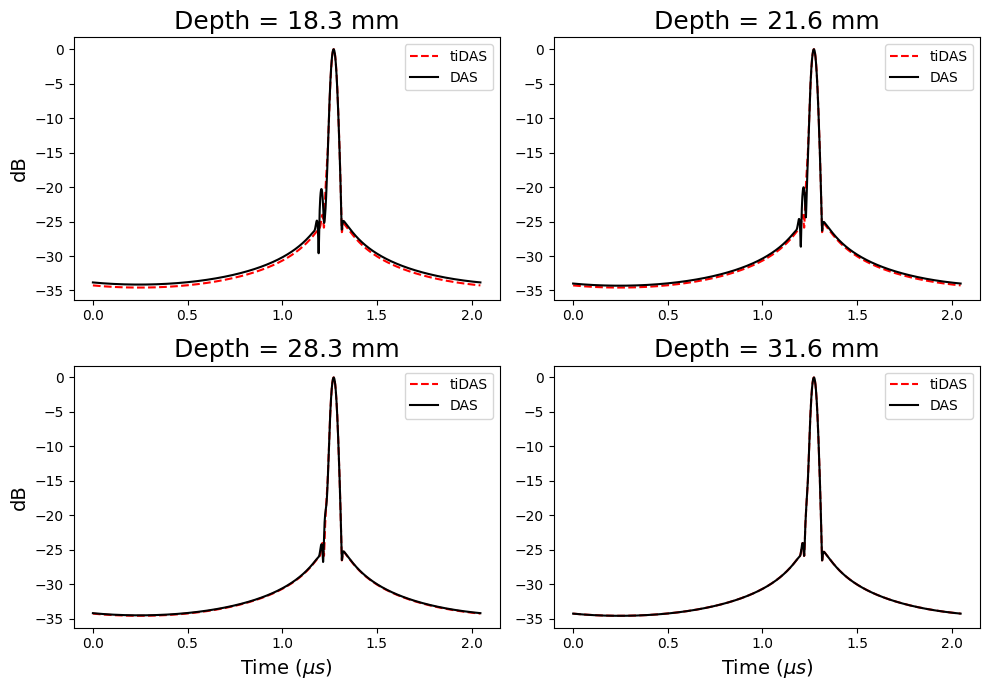}
    \caption{Reconstructed signals at four depths: DAS (black plain line) and tiDAS (red dashed line) show high similarity across all depths. The time axis is relative to the reference depth displayed.}
    \label{fig:PSF}
\end{figure}
Figure~\ref{fig:PSF} presents a qualitative comparison between DAS and tiDAS reconstructions at four representative depths. The results show that both methods yield visually indistinguishable profiles in the vicinity of the signal peak, indicating that the time-invariant approximation effectively preserves the key characteristics of the beamformed signal, even across varying imaging depths.

To further explore the link between the delay profiles and their corresponding scaling parameters, we estimate the set of parameters 
$\Theta$ for each fixed delay curve using Equation~\eqref{eqn:theta}. This procedure produces a matrix of scaling weights, where each row corresponds to a particular delay profile (i.e., a fixed reception focus point), and each column represents a reconstruction location along the imaging line.

\begin{figure}
    \centering
    \includegraphics[width=\linewidth]{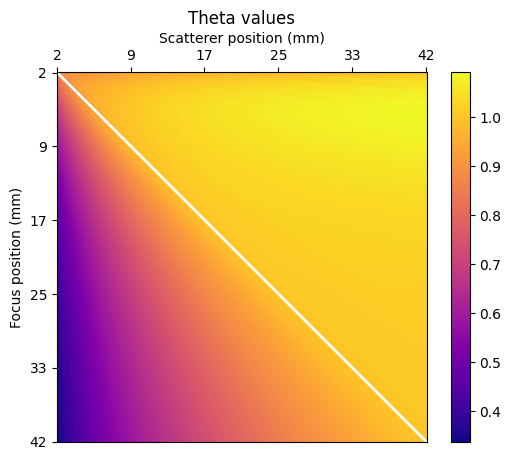}
    \caption{Estimated scaling parameter matrix $\Theta$. Each row corresponds to a fixed delay profile, illustrating its spatial variation.}
    \label{fig:theta}
\end{figure}

\begin{figure}
    \centering
    \includegraphics[width=\linewidth]{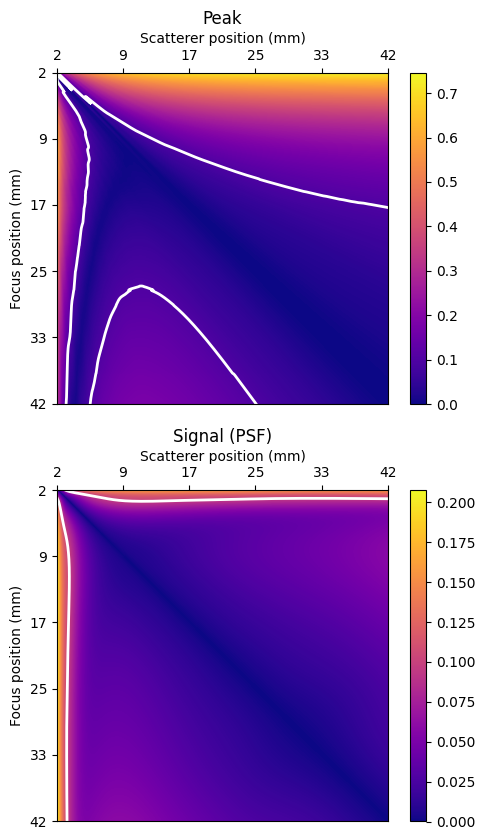}
    \caption{Error analysis: (top) local error (squared $\ell^2$ norm near the peak) and (bottom) global error over the full signal. The blue band marks regions with error below $5\%$.}
    \label{fig:errors}
\end{figure}

The resulting structure, depicted in Figure~\ref{fig:theta}, highlights the smooth spatial variation of the scaling parameters across different configurations. This smoothness supports the validity and practicality of the convolutional tiDAS model, as it suggests that a small number of representative delay profiles and scaling curves may be sufficient to generalize across the imaging domain.
A quantitative evaluation of the reconstruction error is shown in Figure~\ref{fig:errors}. The left panel reports the localized error, computed as the squared $\ell^2$ of the difference between the DAS and tiDAS signals, restricted to a temporal window centered around the peak. This local region is the most relevant for interpreting point scatterers, as it captures the dominant features of the signal.
In the right panel, we present the global reconstruction error, computed across the whole signal duration. Each point in both plots corresponds to a specific combination of scatterer position and fixed delay curve.
To facilitate interpretation, the blue contours in each plot mark the region where the reconstruction error remains below 5$\%$, serving as a practical threshold for acceptable approximation quality. These results demonstrate that tiDAS maintains high fidelity in key regions of interest, while also offering insight into the trade-offs involved in delay curve selection.
These results demonstrate that, for a substantial number of fixed delay curves, it is possible to accurately reconstruct the whole range of scatterer positions with negligible error. This finding underscores the flexibility and efficiency of the tiDAS model, suggesting that a small, strategically selected set of fixed delay profiles can achieve high-quality reconstructions across a broad spatial domain. 

\subsection{Line reconstruction}
To further assess the performance of the time-invariant approximation in a more realistic imaging scenario, we consider the reconstruction of a whole image line containing multiple point scatterers, following the model described in \eqref{eqn:line_conv}. 
Unless otherwise specified, all scatterers are assumed to have unit amplitude, which simplifies the comparison by isolating the effects of the reconstruction method from amplitude scaling.
\begin{figure}
    \centering
    \includegraphics[width=\linewidth]{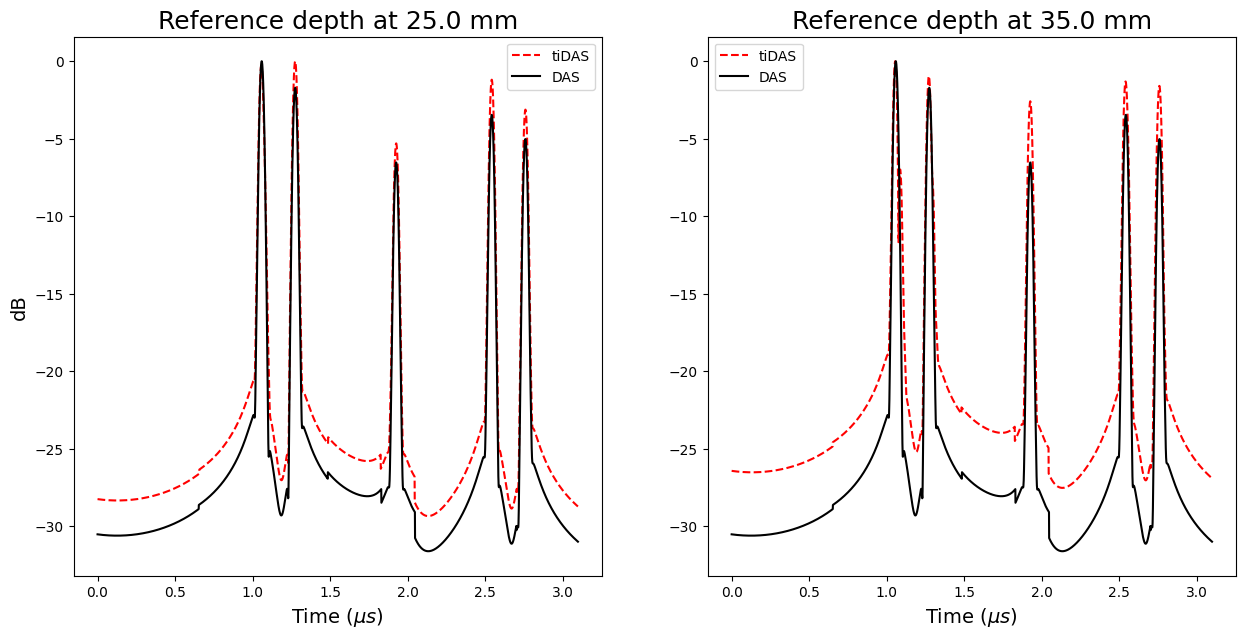}
    \caption{Line reconstruction with five scatterers using two different delay curve at $25$ and $35\,mm$ depth.}
    \label{fig:line}
\end{figure}
The reconstruction is carried out using the delay curve corresponding to a depth of $25\, mm$, a region previously identified as exhibiting low approximation error. 
As illustrated in Figure~\ref{fig:line}, the tiDAS method produces higher peak amplitudes than the standard DAS, indicating enhanced focusing at the scatterer locations. 
Quantitatively, Table~\ref{tab:sidelobes} shows that the average sidelobe (SL) levels for tiDAS are about $2$ to $4, dB$ higher than DAS in both cases of the uniform line with the two different reference depths, with values ranging from about $-26\, dB$ to $24\, dB$ for tiDAS compared to $-28.8\, dB$ for DAS.
Despite this increase, the relative side lobe error is below $0.15$, indicating that the side lobes remain controlled and do not compromise image quality.

\begin{figure}
    \centering
    \includegraphics[width=\linewidth]{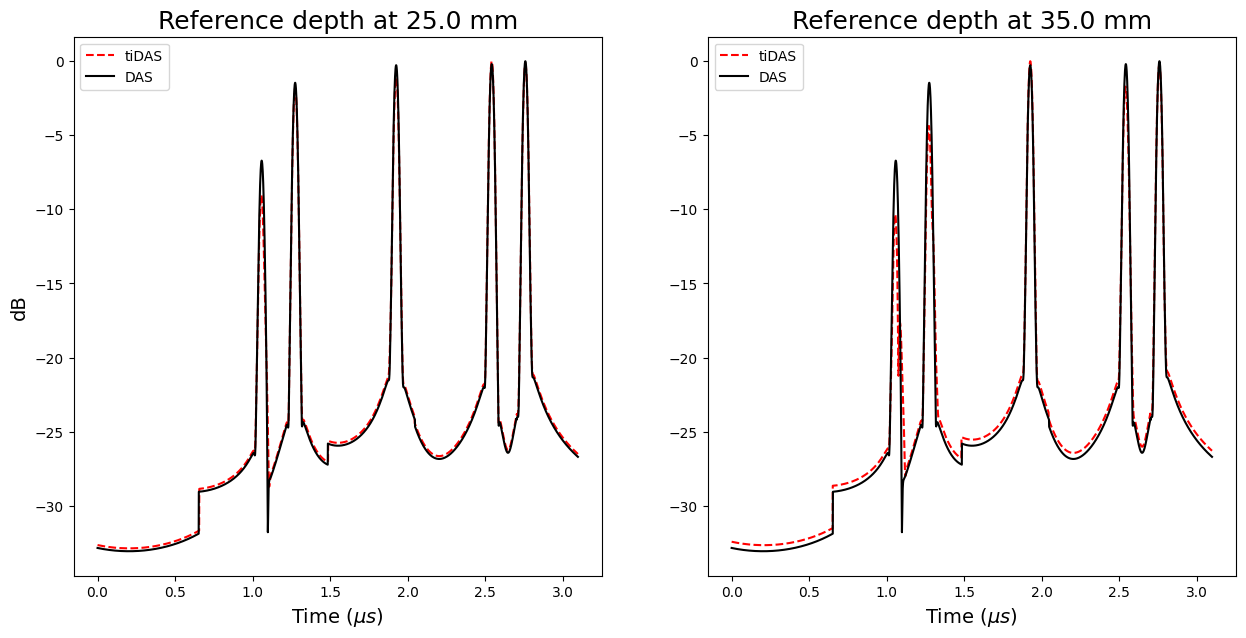}
    \caption{Reconstruction of a line with varying scatterer intensities, increasing with depth. The recontructions are performed using two different delay curve at $25$ and $35\,mm$ depth.}
    \label{fig:intensity}
\end{figure}

Moreover, the increase in side lobe levels is offset by the consistent improvement or maintenance of the main lobe amplitudes observed in tiDAS, resulting in sharper peaks and improved signal contrast. 
This balance ensures that the overall reconstruction quality is not degraded, and in fact benefits from tiDAS is more focused response.
To investigate performance under a more realistic intensity distribution, we also test a configuration where scatterer amplitudes increase with depth, mimicking the natural variability found in biological tissues. 
In this case, shown in Figure~\ref{fig:intensity}, the agreement between tiDAS and DAS futher improves: relative differences decrease, and both peak and side lobe amplitudes align more closely. 
These results demonstrate the robustness of the tiDAS approximation, even in the presence of spatially varying reflectivity.

\begin{figure}
    \centering
    \includegraphics[width=\linewidth]{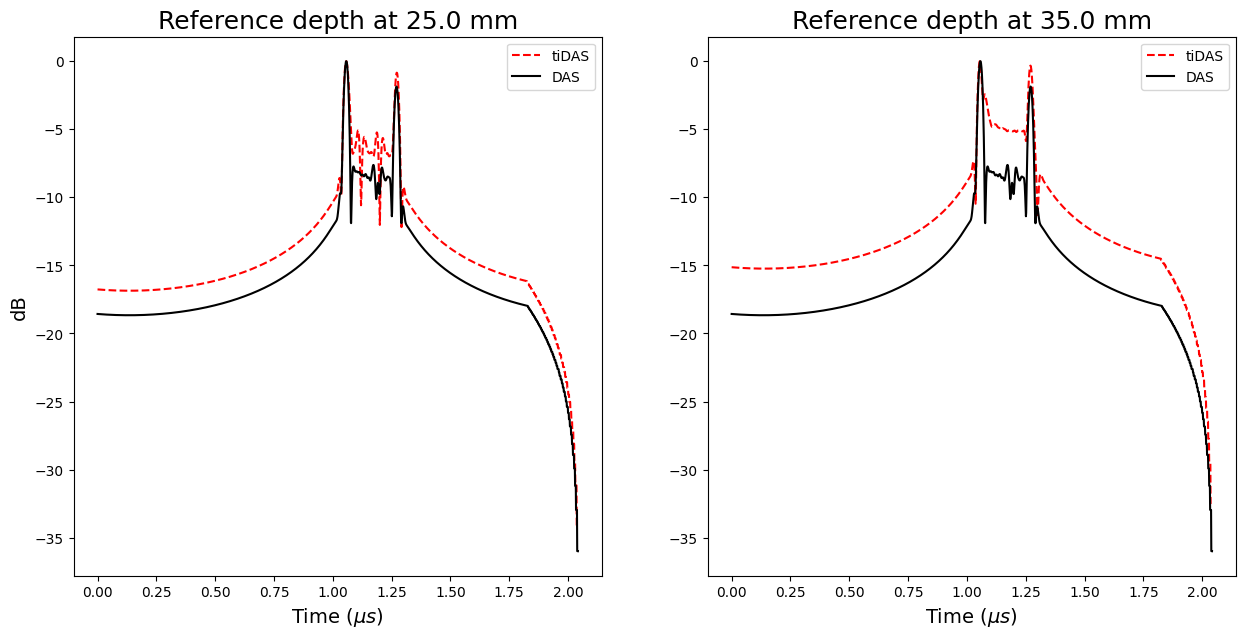}
    \caption{Reconstruction of a line with 100 scatterers using two different delay curve at $25$ and $35\,mm$ depth.}
    \label{fig:line100}
\end{figure}
To further support these observations, we evaluate the reconstruction of a line containing $100$ equally intense scatterers uniformly distributed along the imaging axis. 
As shown in Figure~\ref{fig:line100}, the tiDAS method continues to produce sharper peaks than standard DAS, while faithfully preserving the overall structure and intensity distribution of the signal. 
The mean side lobe level for tiDAS in this denser configuration is approximately $16.07, dB$ (depth $25, mm$), which is about $1.4, dB$ higher than DAS, with a relative side lobe error of $0.10$ (Table~\ref{tab:sidelobes}). 
Even at $35, mm$ depth, the relative error remains moderate at $0.17$. 
This slight increase in side lobe energy, manifesting as a raised baseline, does not compromise the visibility or localization of the main scatterers, confirming that tiDAS maintains reconstruction fidelity in cluttered environments.

These results confirm that the tiDAS approximation remains robust even in densely populated and cluttered environments, reinforcing its suitability for real-time ultrasound imaging in anatomically complex regions.

\begin{table}
\centering
\begin{tabular}{ccccc}
\textbf{Scatterers} & \multicolumn{1}{c}{\textbf{\begin{tabular}[c]{@{}c@{}} Depth \\ (mm)\end{tabular}}}
 & \multicolumn{1}{c}{\textbf{\begin{tabular}[c]{@{}c@{}}Mean DAS\\ SL level\end{tabular}}} & \multicolumn{1}{c}{\textbf{\begin{tabular}[c]{@{}c@{}}Mean tiDAS\\ SL level\end{tabular}}} & \multicolumn{1}{c}{\textbf{\begin{tabular}[c]{@{}c@{}}Relative\\ SL error\end{tabular}}} \\[4pt]
\toprule
5 uniform & 25.0  & $-28.42 \, dB$ & $-26.23 \, dB$ & $0.08$\\[1pt]\midrule
5 uniform  & 35.0  & $-28.41 \, dB$ & $-24.34 \, dB$ & $0.15$\\[1pt]\midrule
5 different  & 25.0 & $-28.83 \, dB$ & $-26.29 \, dB$ & $0.09$\\[1pt]\midrule
5 different  & 35.0 & $-28.88 \, dB$ & $-26.31 \, dB$ & $0.14$\\[1pt]\midrule
100 uniform  & 25.0    & $-17.43 \, dB$ & $-16.07 \, dB$ & $0.10$\\[1pt]\midrule
100 uniform  & 35.0   & $-17.23\, dB$ & $-14.94\, dB$ & $0.17$\\[1pt]\bottomrule
\\
\end{tabular}
\caption{Comparison of mean side lobe (SL) levels between DAS and tiDAS for different line reconstruction scenarios, along with the relative SL error.}
\label{tab:sidelobes}
\end{table}

\subsection{Computational time analysis}
To assess the practical benefits of the proposed time-invariant approximation, we measured the computational times required to perform the four main simulation experiments presented in this study. 
As shown in Table~\ref{tab:timings}, the tiDAS method consistently outperforms the classical DAS algorithm across all configurations. 

This advantage is particularly evident in the reconstruction of linear profiles, where the convolutional structure of tiDAS provides substantial speed-ups.
For example, in the reconstruction of a 100-point line, tiDAS reduces the computational time by nearly eightfold compared to DAS ($0.021 , s$ vs.\ $0.167 , s$), demonstrating its suitability for dense signal reconstructions.

Even in more computationally intensive experiments, such as the simulation involving 600 scatterer positions reconstructed across 600 configurations, tiDAS completes the task in $21.11, min$, compared to $28.35, min$ with DAS.
This highlights that the performance gains scale with problem complexity, reinforcing tiDAS as a promising approach for high-throughput or real-time imaging applications.

\begin{table}
\centering
\begin{tabular}{ccc}
       \textbf{Experiment}  & \textbf{DAS} & \textbf{tiDAS} \\[4pt]
       \toprule
       Four points in one configuration & $0.061 \, s$ & $0.020 \, s$ \\[1pt]\midrule
       600 points in 600 configurations & $28.35 \, min$ & $21.11\, min$ \\[1pt]\midrule
       Line with five uniform scatterers & $0.020 \, s$ & $0.006 \, s$ \\[1pt]\midrule
       Line with five different scatterers & $0.016 \, s$ & $0.005 \, s$ \\[1pt]\midrule
       Line with 100 uniform scatterers & $0.167 \, s$ & $0.021 \, s$ \\[1pt]\bottomrule\\  
\end{tabular}
\caption{Comparison of computational times between DAS and tiDAS for the main experiments.}
\label{tab:timings}
\end{table}
It is important to note that these measurements were obtained in a synthetic simulation environment. 
As such, they do not account for hardware-specific factors that may affect runtime in real-world systems, such as memory access patterns, I/O overhead, and GPU acceleration capabilities. 
Future studies will be needed to validate the practical benefits of tiDAS in clinical imaging pipelines, where these factors could play a significant role.

\section{Conclusion}
In this work, we introduced a time-invariant approximation of the Delay-and-Sum (DAS) beamforming algorithm, termed tiDAS, with the goal of accelerating ultrasound image reconstruction while maintaining high fidelity. By reformulating the DAS operation as a spatially invariant convolution and associating it with a set of learned scaling parameters $\Theta$, we demonstrated that it is possible to significantly reduce computational complexity without compromising reconstruction quality.

Through extensive simulations, we showed that tiDAS closely matches the performance of standard DAS across a wide range of scatterer positions, with both local and global errors remaining below 5$\%$ in many configurations. The parameter estimation procedure revealed delay curves that perform well over broad spatial regions, suggesting the feasibility of using a small set of predefined profiles for fast image reconstruction.

Furthermore, the reconstruction of whole lines of scatterers confirmed tiDAS ability to preserve key image features, making it suitable for real-time clinical applications. Despite these promising results, the study is currently limited to a one-dimensional reconstruction framework. Extending the approach to two-dimensional (2D) and ultimately full B-mode imaging remains an essential next step, as it would allow for a more realistic assessment of the method performance in practical scenarios.

Therefore, future work will focus on generalizing the convolutional formulation to 2D, evaluating its computational advantages, and exploring strategies for efficiently estimating and applying scaling parameters across the whole image plane. Overall, this approach opens the door to efficient, low-latency ultrasound imaging pipelines, offering exciting opportunities for integrating advanced post-processing and deconvolution methods into time-constrained diagnostic workflows.

\section{Acknowledgements}
FB acknowledges the financial support under the National Recovery and Resilience Plan (NRRP), Mission 4, Component 2, Investment 1.1, Call for tender No. 104 published on 2.2.2022 by the Italian Ministry of University and Research (MUR), funded by the European Union – NextGenerationEU– Project Title “Computational mEthods for Medical Imaging (CEMI) – CUP D53D23005830006 - Grant Assignment Decree No. 973 adopted on 30.6.2023 by the Italian Ministry of Ministry of University and Research (MUR).
CR acknowledge the support obtained by INdAM - GNCS Project 'Approcci RBF per l'inpainting di problemi PDE, adattati con reti neurali CUP\_E53C24001950001.
This work was also supported by NextGenerationEU (NGEU) and funded by the Ministry of University and Research (MUR), National Recovery and Resilience Plan (NRRP), project RAISE (ECS00000035)—Robotics and AI for Socio-economic Empowerment (DN. 1053 del 23.06.2022). The views and opinions expressed herein are those of the authors alone and do not necessarily reflect those of the European Union or the European Commission. Neither the European Union nor the European Commission can be held responsible for them. MP also
acknowledges the financial support of the National Recovery and Resilience Plan (NRRP), Mission
4, Component 2, Investment 1.1, Call for tender No. 104 published on 2.2.2022 by the Italian Ministry of University and Research (MUR), funded by the European Union – NextGenerationEU –
Project Title “Inverse Problems in Imaging Sciences (IPIS)” – CUP D53D23005740006 – Grant Assignment Decree No. 973 adopted on 30/06/2023 by the Italian Ministry of Ministry of University
and Research (MUR).

\bibliographystyle{alpha}
\bibliography{Razzetta}

\end{document}